\documentclass[conference]{IEEEtran}
\IEEEoverridecommandlockouts

\usepackage{subfigure}
\usepackage{threeparttable}
\usepackage{booktabs}
\usepackage{multirow}

\usepackage{cite}
\usepackage{amsmath,amssymb,amsfonts}
\usepackage{algorithmic}
\usepackage{graphicx}
\usepackage{textcomp}
\usepackage{xcolor}
\def\BibTeX{{\rm B\kern-.05em{\sc i\kern-.025em b}\kern-.08em
    T\kern-.1667em\lower.7ex\hbox{E}\kern-.125emX}}

\begin{document}

    \title{
    Design of Transit Networks: Global Optimization of Continuous Approximation Models via Geometric Programming \\
    \thanks{This research is supported by the Natural Science Foundation of China (NSFC, Grant No. 725B2037).}
    \thanks{* Haoyang Mao is the corresponding author.}
    }

    \author{\IEEEauthorblockN{1\textsuperscript{st} Haoyang Mao*}
    \IEEEauthorblockA{\textit{Department of Electrical and Electronic Engineering} \\
    \textit{The Hong Kong Polytechnic University}\\
    Hong Kong, PR China \\
    hyang.mao@connect.polyu.hk}
    \and
    \IEEEauthorblockN{2\textsuperscript{nd} Weihua Gu}
    \IEEEauthorblockA{\textit{Department of Electrical and Electronic Engineering} \\
    \textit{The Hong Kong Polytechnic University}\\
    Hong Kong, PR China \\
    weihua.gu@polyu.edu.hk}
    \and
    \IEEEauthorblockN{3\textsuperscript{rd} Wenbo Fan}
    \IEEEauthorblockA{\textit{Department of Electrical and Electronic Engineering} \\
    \textit{The Hong Kong Polytechnic University}\\
    Hong Kong, PR China \\
    wenbo.fan@polyu.edu.hk}
    \and
    \IEEEauthorblockN{4\textsuperscript{th} Zhicheng Jin}
    \IEEEauthorblockA{\textit{Department of Electrical and Electronic Engineering} \\
    \textit{The Hong Kong Polytechnic University}\\
    Hong Kong, PR China \\
    zhicheng.jin@connect.polyu.hk}
    \and
    \IEEEauthorblockN{5\textsuperscript{th} Xiaokuan Zhao}
    \IEEEauthorblockA{\textit{Department of Electrical and Electronic Engineering} \\
    \textit{The Hong Kong Polytechnic University}\\
    Hong Kong, PR China \\
    xiaokuan98.zhao@connect.polyu.hk}
    }

	\maketitle
	\begin{abstract}
		Continuous approximation (CA) models have been widely adopted in transit network design studies due to their strong analytical tractability and high computational efficiency. However, such models are typically formulated as nonconvex optimization problems, and existing solution approaches mainly rely on iterative algorithms that exploit first-order optimality information or nonlinear programming solvers, whose solution quality lacks stability guarantees under complex demand conditions. This paper proposes a geometric programming (GP)–based CA method for transit network design
        , which can be efficiently solved to global optimality. Numerical experiments are conducted on both homogeneous and heterogeneous network settings to evaluate the effectiveness of the proposed approach. Comprehensive tests are performed under the combinations of six heterogeneous demand distributions, four levels of total passenger demand, and three value-of-time parameters. The results indicate that the GP approach consistently outperforms the coordinate descent method across all test cases, achieving cost reductions of approximately 1\%–4\%, even when the latter converges to identical solutions under different initializations. In comparison, nonlinear programming solvers, with fmincon as a representative example, are able to obtain globally optimal solutions comparable to those of the GP approach in low-demand heterogeneous networks; however, their performance becomes unstable under high-demand conditions. These findings demonstrate that GP provides an efficient and robust optimization framework for solving CA–based transit network design problems, especially in high-demand and highly heterogeneous network environments.
	\end{abstract}

	\section{Introduction}
	\label{sec1}
    Transit networks constitute the backbone of urban public transportation systems and play a critical role in both residents' travel efficiency and the operational costs of public transit services. With the increasing complexity of travel demand patterns, traditional optimization approaches based on discrete network flow models have revealed inherent limitations when addressing large-scale transit network design problems, most notably the NP-hard and curse of dimensionality. These methods often struggle to strike an effective balance between computational tractability and solution quality. Against this backdrop, the continuous approximation (CA) approach has gradually emerged as an important theoretical framework for transit network design and optimization due to its strong analytical tractability and high computational efficiency.

    The fundamental idea of CA is to transform discrete elements—such as lines, stops, and origin–destination (OD) demand matrices—into continuous representations of density per unit area, expressed as functions of time and space \cite{holroydOPTIMUMBUSSERVICE1967,chienOptimizationGridTransit1997}. This transformation not only substantially reduces computational complexity but also enables a systematic derivation of analytical relationships between network design parameters and system performance measures. As a result, CA-based models provide transit planners with theoretical insights that are both interpretable and policy-relevant.
    

    However, the frameworks of these studies are typically formulated as nonlinear programming (NLP) problems with non-convex objective functions and constraints. Existing solution approaches largely rely on iterative algorithms that exploit first-order optimality information, such as gradient descent and coordinate descent methods \cite{chienOptimizationGridTransit1997,ouyangContinuumApproximationApproach2014,chenOptimalTransitService2015}, or general-purpose nonlinear optimization solvers (e.g., fmincon) \cite{wuOptimalDesignTransit2020b,zhenComparingSkipstopAllstop2025}. These methods suffer from shortcomings. On the one hand, gradient-based algorithms, while offering numerical stability, are highly susceptible to local optima and thus cannot guarantee solution quality. On the other hand, general NLP solvers often exhibit pronounced instability, with solution quality being highly sensitive to the choice of initial points. The difficulty of obtaining reliable globally optimal solutions has, to some extent, hindered the broader practical application of CA methods in transit network design.
    
    To overcome these challenges, this paper introduces geometric programming (GP) as an optimization paradigm, which is a specialized optimization technique tailored for problems with posynomial structures. Its core theoretical foundation lies in the fact that, through exponential variable substitutions and logarithmic transformations, a non-convex problem can be equivalently converted into a convex optimization problem. This transformation enables the use of efficient interior-point methods to obtain globally optimal solutions in polynomial time, without requiring any initial solution guesses \cite{avrielAdvancesGeometricProgramming1980}. Importantly, this process preserves the physical interpretation of the original variables while providing strong theoretical guarantees of global optimality.
    
    GP has been successfully applied in various engineering fields, including communication networks \cite{chiangPowerControlGeometric2007}, integrated circuit design \cite{boydDigitalCircuitOptimization2005}, and natural gas infrastructure networks \cite{misraOptimalCompressionNatural2015}, demonstrating its unique advantages in handling large-scale optimization problems with multiple constraints. However, to the best of the authors' knowledge, the application of GP to CA models for transit network design has not yet been reported in the literature. This methodological gap constitutes the primary motivation of the present study. In summary, this study makes the following contributions:
    
    \begin{itemize}
    \item The classical CA models for transit network design are solved through GP, thereby providing a solution paradigm with guaranteed global optimality.
    \item Extensive numerical experiments are conducted under both homogeneous and heterogeneous networks, and the proposed GP-based approach is compared with gradient-based algorithms and commonly used nonlinear programming solvers.
    \end{itemize}

    The paper unfolds as follows: Section 2 presents the CA modeling framework for homogeneous and heterogeneous network design problems, and section 3 shows the solution approach via GP. Section 4 describes the numerical case studies. Section 5 summarizes the findings and outlines directions for future research.
	
	\section{Optimal Network Design Model}
	\label{sec2}
\begin{figure}[htbp]
    \centering
    \subfigure[Homogeneous network]{
        \includegraphics[width=0.45\linewidth]{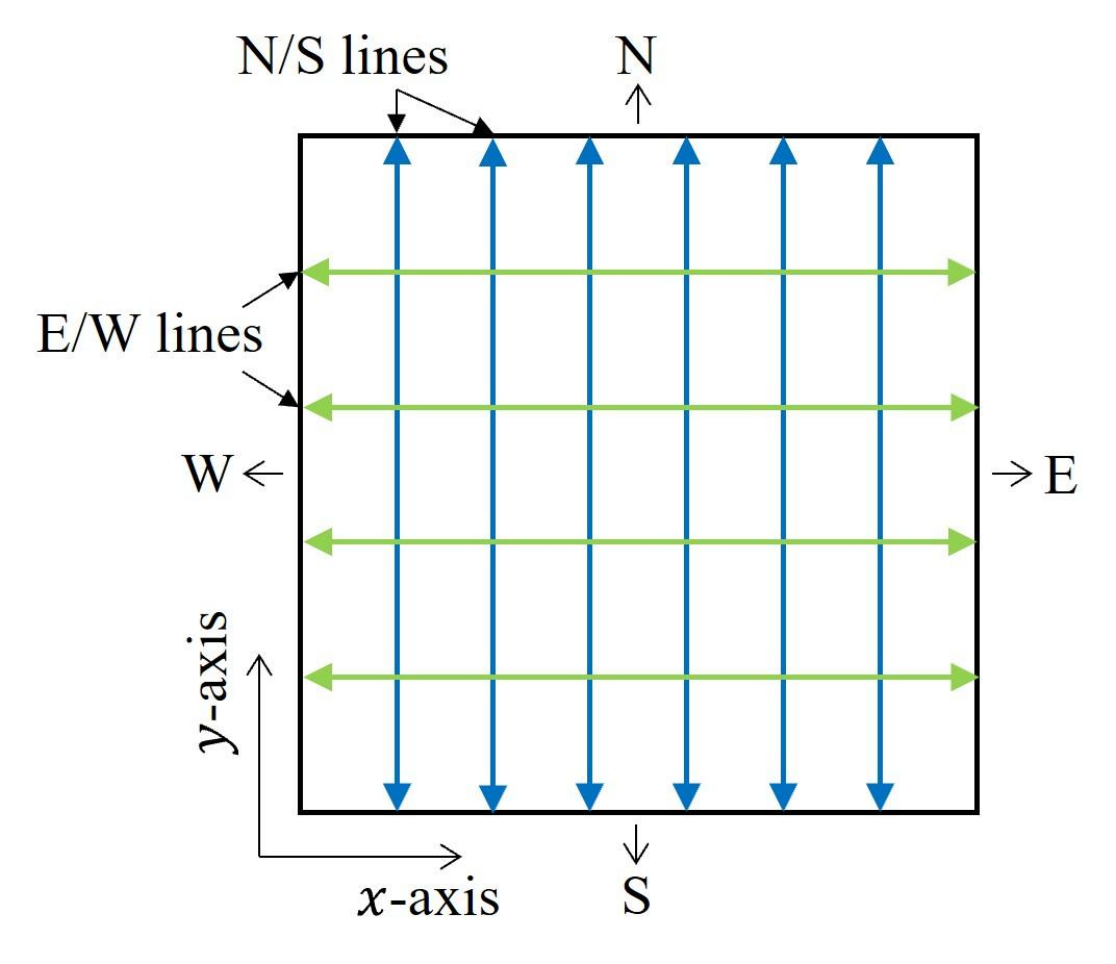}
        \label{Fig_Network_Homnet}
    }
    \subfigure[Heterogeneous network]{
        \includegraphics[width=0.45\linewidth]{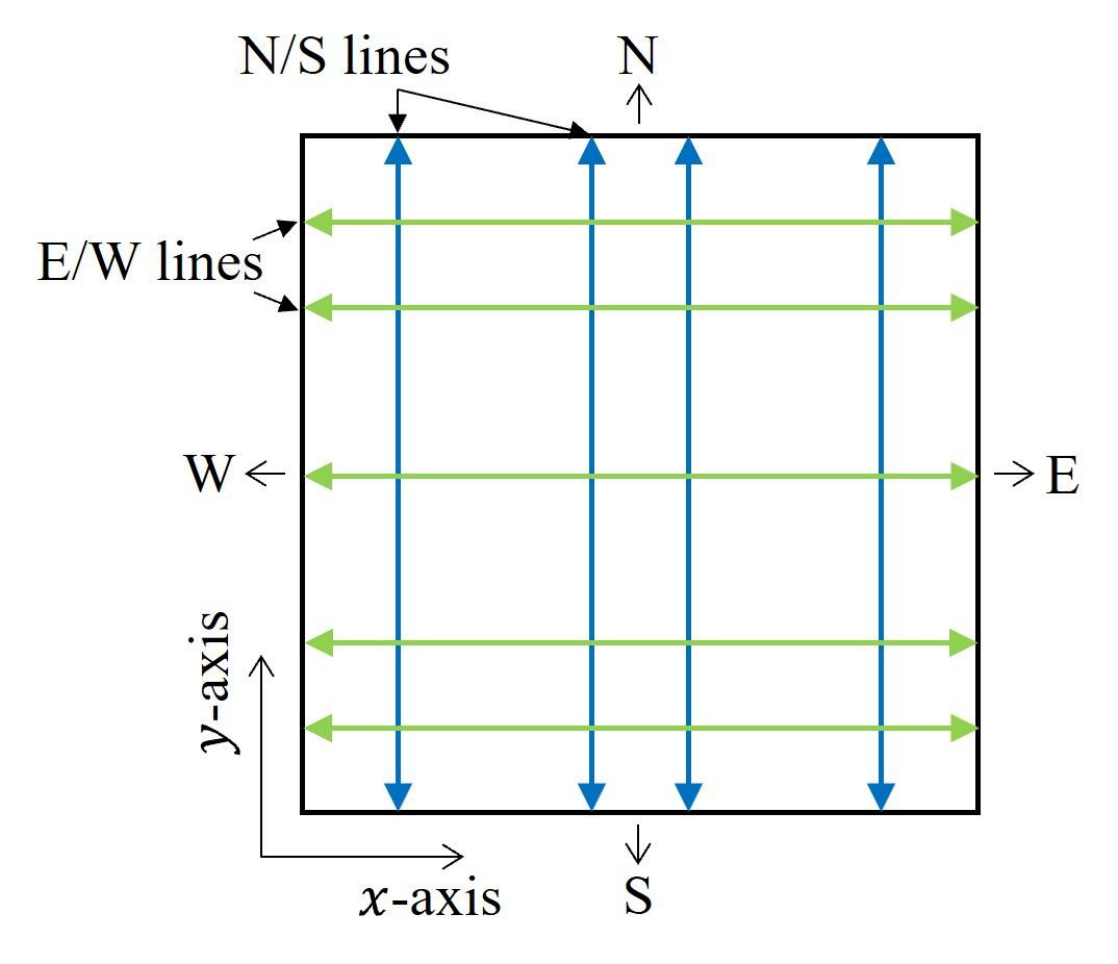}
        \label{Fig_Network_Hetnet}
    }
    \caption{Illustration of homogeneous and heterogeneous grid-based transit networks in a square city.}
    \label{Fig_Network}
    \end{figure}

    This study considers a square city as the study area and constructs both homogeneous and heterogeneous grid-based transit networks, as illustrated in Fig. \ref{Fig_Network}. The study area is denoted by $R$, with a side length of $\left| R \right|$. A location within $R$ is represented by the two-dimensional vector $\left( x,y \right)$, where $x,y \in \left[ 0,\left| R \right|\right]$. Model variables and parameters are labeled by operational direction $i \in I$, where set $I$ includes $E$ (Eastbound), $W$ (Westbound), $N$ (Northbound), and $S$ (Southbound), respectively.

    In a homogeneous network, the line density and service headway in each operational direction are assumed to be constant over the entire study area $R$.
    In contrast, the heterogeneous network considered in this study, following Chien and Schonfeld (1997) \cite{chienOptimizationGridTransit1997}, allows the line density and service headway in each direction $i \in I=\left\{ E,W,N,S \right\}$ to vary spatially in the direction perpendicular to the line orientation, while remaining invariant along the line direction
    Moreover, no lateral detours in the perpendicular direction are permitted in the heterogeneous network.

    The model parameters are distinguished between spatially varying (localized) quantities and spatially invariant parameters. To establish a unified modeling framework applicable to both homogeneous and heterogeneous networks, the localized parameters are formulated as functions of the spatial coordinates $\left( x,y \right)$. All parameters are summarized in Appendix \ref{app1}. The model is subsequently simplified in Section \ref{Sec2_4} according to the characteristics of each network type.

	\subsection{Basic concepts and assumptions}
    \label{sec2_1}
    The model involves two categories of decision variables: the line density $\delta_i\left(x,y\right)$ and the service headway $h_i\left(x,y\right)$. The former determines the structural configuration of the transit network, while the latter governs the operational characteristics of transit vehicles. Reflecting real-world practice, most transit lines operate in bidirectional symmetry. Accordingly, we assume $\delta_E\left(x,y\right)=\delta_W\left(x,y\right),\delta_N\left(x,y\right)=\delta_S\left(x,y\right)$; and $h_E\left(x,y\right)=h_W\left(x,y\right), h_N\left(x,y\right)=h_S\left(x,y\right)$, so as to reduce operational complexity.

    Regarding stop placement, all stops are assumed to be located at the intersections of any two perpendicular transit lines. Passengers are allowed to transfer between orthogonal lines at these stops, and vehicles are assumed to operate without skip-stop service. Intermediate stops between adjacent transfer stops are neglected in this study and are left for future investigation.
    
    Based on the two categories of decision variables, the vehicle flow of lines in direction $i$ per unit cross-sectional distance in the network, denoted by $q_i\left(x,y\right)$, can be computed as shown in Eq. (\ref{Eq_Vehicle_Flow}):
    
    \begin{eqnarray}\label{Eq_Vehicle_Flow}
    q_i\left(x,y\right) = \frac{\delta_i\left(x,y\right)}{h_i\left(x,y\right)}, i \in I, x,y \in \left[0,\left|R\right|\right]
    \end{eqnarray}

    Passenger demand density is treated as an exogenous input to the model and is denoted by $\lambda_i\left(x_o,y_o,x_d,y_d\right)$, which is a function of four spatial coordinates. Here, $\left(x_o,y_o\right)$ and $\left(x_d,y_d\right)$ represent the origin and destination of the passenger, respectively. With respect to their travel behavior, following previous studies \cite{chriquiCommonBusLines1975,daganzoPublicTransportationSystems2019,fanOptimalDesignIntersecting2018a,zhichengjinMultiobjectiveOptimizationModel2024}, the following assumptions are adopted:
    
    \begin{itemize}
    \item Passengers walk to the nearest transit stops for boarding and alighting.
    \item Passengers arrive at the origin stop randomly, without consulting the transit schedule.
    \item Passenger service follows the first-in, first-served (FIFS) principle, since passenger flows do not exceed the service capacity provided by the vehicle flows in any direction $i$ at any location $\left(x,y\right)$ in the network (see the constraint formulation in Section \ref{Sec2_3}).
    \end{itemize}

    Under these assumptions, the high-dimensional passenger demand density function can be aggregated to construct several categories of demand functions, including the boarding passenger density $ \lambda_{bo}^i\left(x,y \right)$, alighting passenger density $\lambda_{al}^i\left(x,y \right)$, on-board flux density $\lambda_{fl}^i\left(x,y \right)$, and transferring density $\lambda_{tr}^i\left(x,y \right)$. The detailed derivation of these quantities follows the methodology proposed by Chien and Schonfeld (1997) \cite{chienOptimizationGridTransit1997}.
    
    \subsection{Network design objective function}
    \label{Sec2_2}    
    The total system cost of the transit network is denoted by $Z$, which consists of two components: the agency cost ($Z_A$) and the passenger cost ($Z_P$). The computation of the system cost is given in Eq. (\ref{Eq_Z}), where $\mu$ represents passengers' average value of time (VOT).

    \begin{eqnarray}\label{Eq_Z}
    Z = \frac{1}{\mu}Z_A + Z_P
    \end{eqnarray}

    Agency cost, $Z_A$, comprises four metrics: (i) the line infrastructure length, $N_l$; (ii) the number of stops, $N_s$; (iii) the vehicle distance traveled, $N_k$; and (iv) the vehicle time traveled, $N_h$. Following Daganzo and Ouyang (2019) \cite{daganzoPublicTransportationSystems2019}, $Z_A$ is computed as shown in Eq. (\ref{Eq_Z_A}), where the parameters $\pi_l,\pi_s,\pi_k$, and $\pi_h$ denote unit costs associated with the four metrics, $\bar{i}$ represents the direction perpendicular to direction $i$, $v$ denotes cruising speed of transit vehicles, and $\tau$ is the fixed delay per stop.

    \begin{subequations}\label{Eq_Z_A}
    \begin{align}
    Z_A &= \pi_l N_l + \pi_s N_s + \pi_k N_k + \pi_h N_h , \label{Eq_Z_A_total}\\
    N_l &= \sum_{i \in I} \int_{y=0}^{\left|R\right|} \int_{x=0}^{\left|R\right|} \delta_i\left(x,y\right) dx \, dy , \label{Eq_Z_A_N_l}\\
    N_s &= \sum_{i \in I} \int_{y=0}^{\left|R\right|} \int_{x=0}^{\left|R\right|} \delta_i\left(x,y\right) \delta_{\bar{i}}\left(x,y\right) dx \, dy , \label{Eq_Z_A_N_s}\\
    N_k &= \sum_{i \in I} \int_{y=0}^{\left|R\right|} \int_{x=0}^{\left|R\right|} q_i\left(x,y\right) dx \, dy , \label{Eq_Z_A_N_k}\\
    N_h &= \sum_{i \in I} \int_{y=0}^{\left|R\right|} \int_{x=0}^{\left|R\right|} q_i\left(x,y\right) \left( \frac{1}{v}+\tau \delta_{\bar{i}}\left(x,y\right) \right) dx \, dy . \label{Eq_Z_A_N_h}
    \end{align}
    \end{subequations}

    Equations (\ref{Eq_Z_A_N_l}) and (\ref{Eq_Z_A_N_k}) are straightforward, as they can be computed by integrating the corresponding density functions over the study area $R$. As discussed in Section \ref{sec2_1}, stops are assumed to be located at the intersections of any two perpendicular transit lines. Accordingly, the term $\delta_i\left(x,y\right) \delta_{\bar{i}}\left(x,y\right)$ in Eq. (\ref{Eq_Z_A_N_s}) provides an estimate of the stop density at location $\left(x,y\right)$, and the resulting stop density is identical for all four operational directions. Equation (\ref{Eq_Z_A_N_h}) consists of cruising time and dwell time at stops.

    Passenger cost, $Z_P$, is characterized by four primary metrics: (i) the access/egress time, $T_a$; (ii) the out-of-vehicle wait time, $T_w$; (iii) the in-vehicle ride time, $T_r$; and (iv) the transfer penalty, $T_t$. It should be noted that fare payments are excluded from the system cost, as they represent transfers between passengers and the transit agency and therefore cancel out at the system level. $Z_P$ is formulated in Eq. (\ref{Eq_Z_P}), where $v_w$ denotes the walking speed of passengers, and $\sigma$ represents the transfer penalty. Additionally, walking time is multiplied by a perceived walking time factor $\beta_w$, reflecting the higher perceived cost of walking relative to in-vehicle travel \cite{wardmanReviewBritishEvidence2001}.

    \begin{subequations}\label{Eq_Z_P}
    \begin{align}
    Z_P &= T_a + T_w + T_r + T_t ,\label{Eq_Z_P_total}\\
    \nonumber
    T_a &= \beta_w \sum_{i \in I} \int_{0}^{|R|} \int_{0}^{|R|} \left(\lambda_{bo}^i\left(x,y\right)+\lambda_{al}^i\left(x,y\right)\right) \\
    &\qquad \times \left(
    \frac{1}{4 v_w \delta_i\left(x,y\right)} +
    \frac{1}{4 v_w \delta_{\bar{i}}\left(x,y\right)}
    \right) dx\,dy,
    \label{Eq_Z_P_T_a}\\
    T_w &= \sum_{i \in I} \int_{y=0}^{\left|R\right|} \int_{x=0}^{\left|R\right|} \left(\lambda_{bo}^i\left(x,y\right)+\lambda_{tr}^i\left(x,y\right)\right) \frac{h_i\left(x,y\right)}{2} dx \, dy , \label{Eq_Z_P_T_w}\\
    T_r &= \sum_{i \in I} \int_{y=0}^{\left|R\right|} \int_{x=0}^{\left|R\right|} \lambda_{fl}^i\left(x,y\right) \left( \frac{1}{v}+\tau \delta_{\bar{i}}\left(x,y\right) \right) dx \, dy , \label{Eq_Z_P_T_r}\\
    T_t &= \sigma \sum_{i \in I} \int_{y=0}^{\left|R\right|} \int_{x=0}^{\left|R\right|} \lambda_{tr}^i\left(x,y\right) dx \, dy . \label{Eq_Z_P_T_t}
    \end{align}
    \end{subequations}
    
    Equation (\ref{Eq_Z_P_T_a}) estimates passengers' access and egress time using numerical integration, where $\frac{1}{\delta_i\left(x,y\right)}$ represents the spacing between transit lines in direction $i$ at location $\left(x,y\right)$. Equation (\ref{Eq_Z_P_T_w}) accounts for the total waiting time, including both initial waiting and transfer waiting. Equation (\ref{Eq_Z_P_T_r}) captures passengers' in-vehicle travel time, consisting of cruising time and dwell time at stops. Finally, Equation (\ref{Eq_Z_P_T_t}) directly integrates the transfer density and multiplies it by the corresponding penalty time.
    
    \subsection{Model constraints}
    \label{Sec2_3}
    The model is subject to two primary sets of constraints. First, in any direction $i$ at any location $\left(x,y\right)$ within the study area $R$, the passenger flow must not exceed the maximum service capacity provided by the corresponding vehicle flow, as expressed in Eq. (\ref{Eq_constrant_capacity}), where $C$ denotes the capacity of a transit vehicle.

    \begin{eqnarray}\label{Eq_constrant_capacity}
    \frac{\lambda_{fl}^i\left(x,y\right)}{q_i\left(x,y\right)} \le C, i \in I, x,y \in \left[0,\left|R\right|\right]
    \end{eqnarray}

    Second, non-negativity constraints are imposed on the decision variables to ensure their physical interpretability. Moreover, since both $\delta_i(x,y)$ and $h_i(x,y)$ appear in the denominator of certain terms in the model formulation, they are further restricted to be strictly positive. The resulting constraints are given by

    \begin{eqnarray}\label{Eq_constrant_non_negative}
    \delta_i\left(x,y\right),h_i\left(x,y\right)>0, i \in I, x,y \in \left[0,\left|R\right|\right]
    \end{eqnarray}
    
    \subsection{Model specialization for heterogeneous and homogeneous networks}
    \label{Sec2_4}
    For heterogeneous networks, the two categories of decision variables in each direction $i \in I=\left\{ E,W,N,S \right\}$ are assumed to remain invariant along the corresponding line direction. In addition, considering the bidirectional symmetry of the transit network, they can be specified as:

    \begin{equation}\label{Eq_specialization_hetnet_delta}
    \delta_i(x,y)=
    \left\{
    \begin{aligned}
    \delta_{EW}(y), &\; i \in \{E,W\}
    \\
    \delta_{NS}(x), &\; i \in \{N,S\}
    \end{aligned}
    \right.
    \end{equation}

    Under this specification, $\delta_{EW}(y)$, and $\delta_{NS}(x)$ are treated as the new decision variables. The same specialization applies to the service headway $h_i(x,y)$.

    The objective function must be integrated along the direction of operation for each $i \in I=\left\{ E,W,N,S \right\}$. Taking Eq. (\ref{Eq_Z_A_N_h}) as an example, the resulting expression after integration is given by Eq. (\ref{Eq_specialization_hetnet_objective}), where $q_{EW}(y) = \frac{\delta_{EW}(y)}{h_{EW}(y)}$, and $q_{NS}(x) = \frac{\delta_{NS}(x)}{h_{NS}(x)}$. The transformations of the remaining objective function terms in Eqs. (\ref{Eq_Z_A}) and (\ref{Eq_Z_P}) follow the same procedure and are therefore omitted for brevity.
    \begin{equation}\label{Eq_specialization_hetnet_objective}
    N_h = 
    \begin{aligned}[t]
    & 2\int_{y=0}^{|R|} q_{EW}(y) \left( \frac{|R|}{v}+\tau \int_{x=0}^{|R|}\delta_{NS}(x)dx \right)dy + \\
    & 2\int_{x=0}^{|R|} q_{NS}(x) \left( \frac{|R|}{v}+\tau \int_{y=0}^{|R|}\delta_{EW}(y)dy \right)dx
    \end{aligned}
    \end{equation}

    Under the above specialization, the capacity constraint in Eq. (\ref{Eq_constrant_capacity}) for directions $E,W$ can be reformulated as Eq. (\ref{Eq_specialization_hetnet_capacity}). The corresponding expressions for directions $N,S$ follow analogously.
    
    \begin{eqnarray}\label{Eq_specialization_hetnet_capacity}
    \left\{
    \begin{aligned}
    \frac{\lambda_{fl}^{EW}(y)}{q_{EW}(y)} & \le C \\
    \lambda_{fl}^{EW}(y) & = 
    \displaystyle
    \max_{\substack{i \in \{E,W\}, x \in [0,|R|]}}
    \lambda_{fl}^i(x,y), y \in [0,|R|]
    \end{aligned}
    \right.
    \end{eqnarray}

    For homogeneous networks, the two categories of decision variables in each operational direction are assumed to be constant over the entire study area $R$. Consequently, $\delta_{EW}, \delta_{NS}, h_{EW}$, and $h_{NS}$ are treated as the new decision variables. The objective function is then integrated over all directions. Taking Eq. (3e) as an example, the corresponding expression is presented in Eq. (\ref{Eq_specialization_homnet_objective}), where $q_{EW} = \frac{\delta_{EW}}{h_{EW}}$, and $q_{NS} = \frac{\delta_{NS}}{h_{NS}}$.
    \begin{equation}\label{Eq_specialization_homnet_objective}
    N_h = |R|^2\left[ q_{EW}\left( \frac{1}{v}+\tau \delta_{NS} \right) + q_{NS}\left( \frac{1}{v}+\tau \delta_{EW} \right) \right]
    \end{equation}

    The capacity constraint in Eq. (\ref{Eq_constrant_capacity}) takes the form of Eq. (\ref{Eq_specialization_homnet_capacity}) for directions $E,W$; an identical form applies to directions $N,S$.
    
    \begin{eqnarray}\label{Eq_specialization_homnet_capacity}
    \left\{
    \begin{aligned}
    &\frac{\lambda_{fl}^{EW}}{q_{EW}} \le C \\
    &\lambda_{fl}^{EW} = 
    \displaystyle
    \max_{\substack{i \in \{E,W\}, x,y \in [0,|R|]}}
    \lambda_{fl}^i(x,y)
    \end{aligned}
    \right.
    \end{eqnarray}

	\section{Solution method via geometric programming}
	\label{Sec3}
    GP is a class of nonlinear and nonconvex optimization problems that possesses a number of attractive theoretical and computational properties. A key advantage of GP lies in the fact that it can be transformed into an equivalent convex optimization problem through a suitable change of variables. As a result, global optimal solutions can be computed efficiently without reliance on initial guesses. Interior-point methods developed for GP are known to admit polynomial-time complexity guarantees \cite{nesterovInteriorPointPolynomialAlgorithms1994}. 

    A brief overview of GP follows. GP admits two equivalent representations: a standard (posynomial) form and a convex form obtained via logarithmic transformation. The standard form is expressed in terms of a special class of functions known as monomials and posynomials. A monomial is defined as a function $g:\mathbb{R}_{++}^n\rightarrow \mathbb{R}$, shown as Eq. (\ref{Eq_gp_monomial_function}), where $d \ge 0$ is a multiplicative constant coefficient, and $a(l)\in \mathbb{R}, l=1,2,...\mathcal{L}$ are real-valued exponents.

    \begin{equation}\label{Eq_gp_monomial_function}
    g(\boldsymbol{r})=d\,r_1^{a^{(1)}}\,r_2^{a^{(2)}}...\,r_\mathcal{L}^{a^{(\mathcal{L})}},  \boldsymbol{r}>0
    \end{equation}

    A posynomial is defined as a finite sum of monomials, presented in Eq. (\ref{Eq_gp_posynomial_function}), where $k = 1, 2,... \mathcal{K}$ is the index of the monomial. For example, expressions such as $\ln{2} \cdot r_1^\pi r_2^{-100}+\frac{r_2}{r_3^4}$ constitute posynomials, whereas expressions involving subtraction (e.g., $r_1-r_2$) do not, due to the nonnegativity requirement on multiplicative coefficients.

    \begin{equation}\label{Eq_gp_posynomial_function}
    f(\boldsymbol{r})=\sum_{k=1}^\mathcal{K} d_k\,r_1^{a_k^{(1)}}\,r_2^{a_k^{(2)}}...\,r_\mathcal{L}^{a_k^{(\mathcal{L})}},  \boldsymbol{r}>0
    \end{equation}

    In its standard form, a GP minimizes a posynomial objective subject to posynomial inequality constraints and monomial equality constraints, all normalized to unity:
    
    \begin{equation}\label{Eq_gp_standard_form}
    \begin{aligned}
    \min_{\boldsymbol{r}} \quad 
    & f_0(\boldsymbol{r}) , \\[0.5ex]
    \text{s.t.} \quad
    & f_u(\boldsymbol{r}) \le 1, \quad u = 1,2,...,\mathcal{U}, \\
    & g_w(\boldsymbol{r}) = 1, \quad w = 1,2,...,\mathcal{W}, \\
    & \boldsymbol{r}>0.
    \end{aligned}
    \end{equation}

    Although this formulation is not convex—since posynomials are generally nonconvex functions—the problem can be equivalently transformed into a convex program by applying a logarithmic change of variables. Specifically, defining $s_l=\ln{r_l},\, b_{uk}=\ln{d_{uk}},\, b_w=\ln{d_w}$ and taking logarithms of the objective and constraints yields the convex form of GP:

    \begin{equation}\label{Eq_gp_convex_form}
    \begin{aligned}
    \min_{\boldsymbol{s}} \quad 
    & \ln{f_0(\boldsymbol{s})} = \ln\sum_{k=1}^{\mathcal{K}_0}e^{\left( \boldsymbol{a}_{0k}^T \boldsymbol{s}+b_{0k} \right)} , \\[0.5ex]
    \text{s.t.} \quad
    & \ln{f_u(\boldsymbol{s})} = \ln\sum_{k=1}^{\mathcal{K}_u}e^{\left( \boldsymbol{a}_{uk}^T \boldsymbol{s}+b_{uk} \right)} \le 0, \quad u = 1,2,...,\mathcal{U}, \\
    & \ln{g_w(\boldsymbol{s})} = \boldsymbol{a}_w^T \boldsymbol{s}+b_w = 0, \quad w = 1,2,...,\mathcal{W} , \\
    & \boldsymbol{s}\in \mathbb{R}.
    \end{aligned}
    \end{equation}
    
    This formulation is a convex optimization problem, as the log-sum-exp function is convex \cite{boydConvexOptimization2004}.

    In the proposed CA-based transit network design problem, the unified objective functions in Eqs. (\ref{Eq_Z_A}) and (\ref{Eq_Z_P}) are explicitly formulated as posynomials of the two categories of decision variables, namely the line density $\delta_i(x,y)$ and the service headway $h_i(x,y)$. For the constant terms that are independent of the decision variables, such as Eq. (\ref{Eq_Z_P_T_t}), they can be equivalently represented as monomials with zero exponents on the decision variables.

    The capacity constraint in Eq. (\ref{Eq_constrant_capacity}) satisfies the requirement that inequality constraints in GP standard form must take a posynomial form. Moreover, Eq. (\ref{Eq_constrant_non_negative}) ensures that all decision variables are strictly positive, which is a fundamental assumption of GP.

    The specialization of the unified model to heterogeneous and homogeneous network settings involves only spatial integration and variable aggregation, and does not alter the posynomial structure of the objective function or the constraints. Consequently, the resulting optimization problems for both network types conform to the standard form of GP. Therefore, the problem is solved using the commercial GP solver MOSEK on a desktop computer with a 2.5 GHz CPU (Intel(R) i5-12500H) and a 16G memory, which efficiently yields the global optimum.

	\section{Numerical case studies}
	\label{Sec4}

	\subsection{Set-up}
	\label{Sec4_1}
    In this section, the proposed models and solution approach are primarily demonstrated using a grid-based city with an area of $10 \times 10\, km^2$. For computational tractability, the city is discretized into uniform square cells of size $\Delta^2$. Each cell is represented by its center point $(x_n,y_{n'})$, where $n,n'=1,2,...,\mathcal{N}=\frac{|R|}{\Delta}$. All density functions and service headway within a cell are assumed to be spatially homogeneous and equal to their values at the corresponding cell center.

    To comprehensively assess the robustness and performance of the proposed method, numerical experiments are conducted for both homogeneous and heterogeneous networks under a wide range of scenarios. Specifically, six heterogeneous demand distributions (monocentric, commute, and four chessboard-type patterns), four levels of total passenger demand (5,000, 10,000, 50,000, and 100,000), and three values of passengers' VOT ($25\frac{\$}{hr}, \, 20\frac{\$}{hr}, \, 5\frac{\$}{hr}$) are considered, resulting in a total of 72 test cases. The generation of heterogeneous demand distributions is detailed in Appendix \ref{app2}. The selected VOT values are intended to represent economic conditions in developed, moderately developed, and developing regions.

    In addition, the globally optimal solutions obtained via GP are compared with those derived from gradient-based algorithms (represented by coordinate descent) and general-purpose nonlinear programming solvers (represented by fmincon). For the coordinate descent method, closed-form solutions for the two categories of decision variables are provided in Appendix \ref{app3}. Other parameter values are taken from Ouyang et al. (2014) \cite{ouyangContinuumApproximationApproach2014}, summarized in Appendix \ref{app1}.
    
	\subsection{Numerical results}
	\label{Sec4_2}
\begin{figure}[htbp]
    \centering
    \subfigure[Coordinate descent method (CD) with $D=5,000$ passengers and $\mu=5 \, \frac{\$}{hr}$]{
        \includegraphics[width=0.9\linewidth]{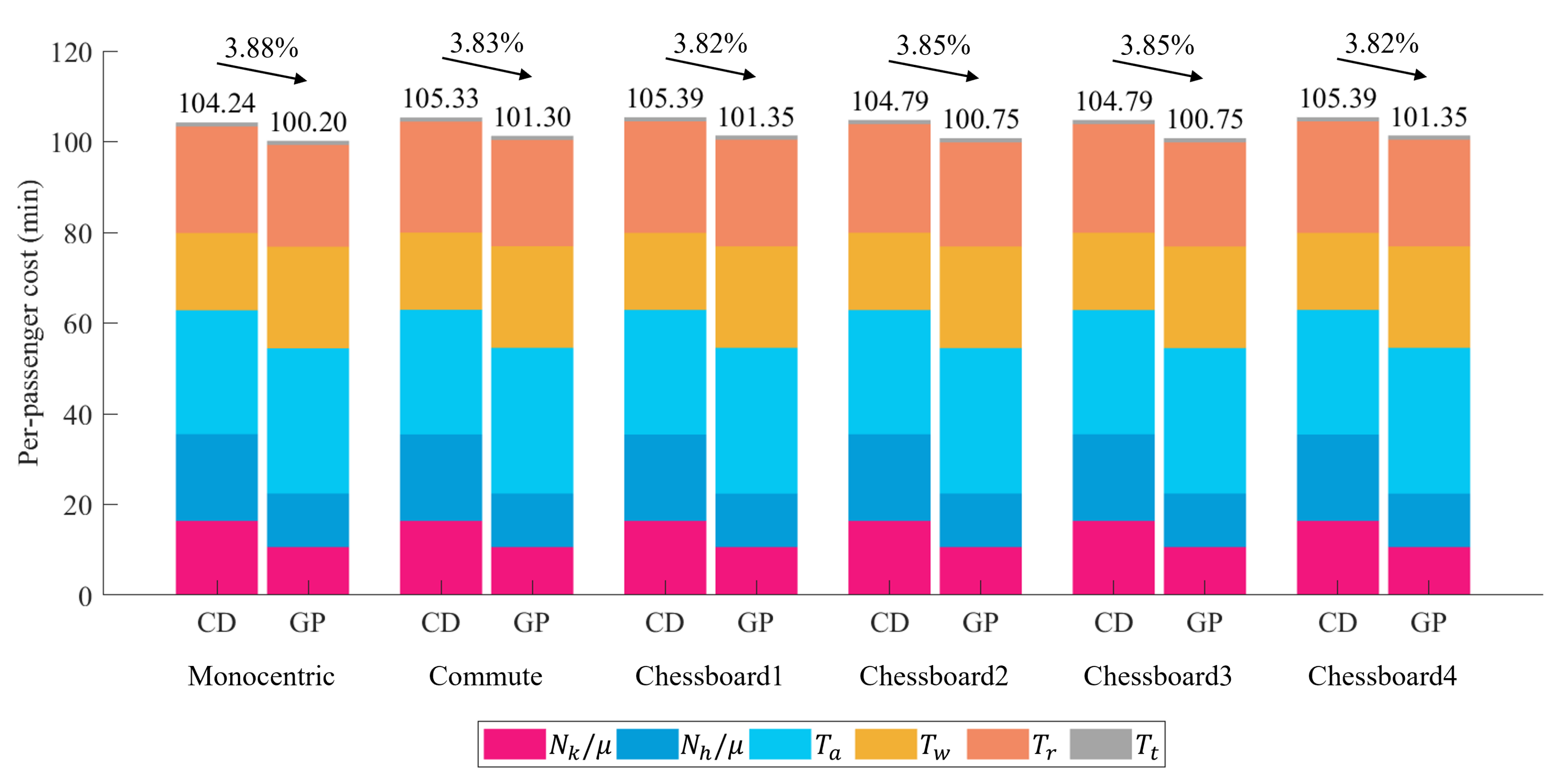}
        \label{Fig_numerical_result_gradient}
    }
    \subfigure[Fmincon with $D=100,000$ passengers and $\mu=5 \, \frac{\$}{hr}$]{
        \includegraphics[width=0.9\linewidth]{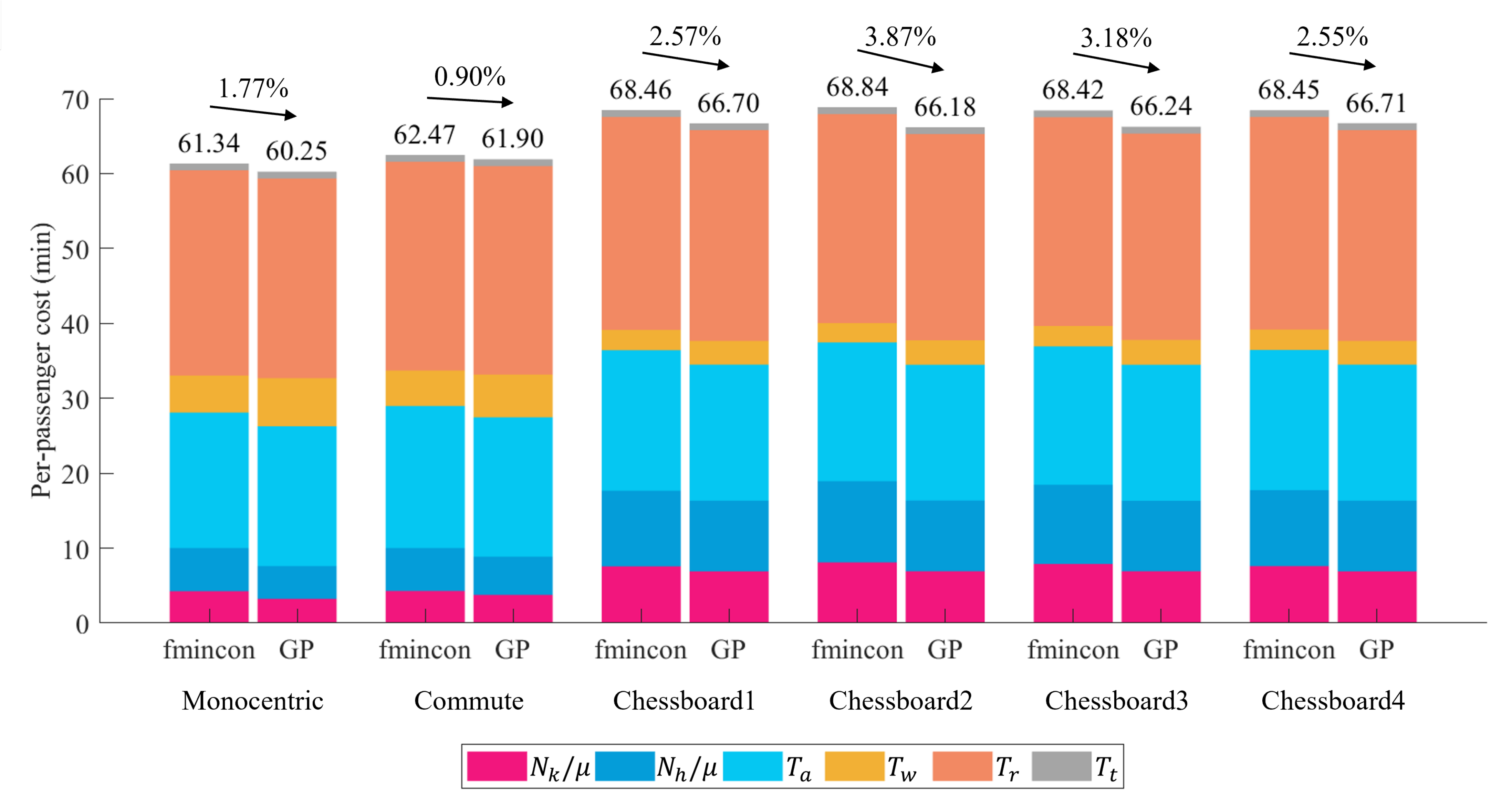}
        \label{Fig_numerical_result_fmincon}
    }
    \caption{Illustration of per-passenger system cost under six demand distributions in heterogeneous network: GP approach vs. two benchmark algorithms.}
    \label{Fig_numerical_result}
    \end{figure}

    Table \ref{Tab_numerical_result} reports the average percentage reductions in total system cost achieved by the GP model, relative to the two benchmark algorithms, across six demand distributions under given total demand levels $D$ and values of $\mu$. Fig. \ref{Fig_numerical_result} illustrates the per-passenger system cost and its individual cost components obtained using the GP approach, in comparison with the two benchmark algorithms.

    \begin{table}[t]
    \caption{Average relative optimality improvement (\%) of the GP approach across six demand distributions for different total demand levels $D$ and values of $\mu$.}
    \label{Tab_numerical_result}
    \centering
    \setlength{\tabcolsep}{3pt}
    \begin{tabular}{c ccc ccc}
    \toprule
    \multirow{3}{*}{$D$} 
    & \multicolumn{6}{c}{Solution method} \\
    \cmidrule(lr){2-7}
    & \multicolumn{3}{c}{Coordinate descent} 
    & \multicolumn{3}{c}{fmincon} \\
    \cmidrule(lr){2-4} \cmidrule(lr){5-7}
    & $\mu=25$ & $\mu=20$ & $\mu=5$ 
    & $\mu=25$ & $\mu=20$ & $\mu=5$ \\
    \midrule
    \multicolumn{7}{c}{\textit{Heterogeneous network}} \\
    \midrule
    5000   & 2.70& 2.86& 3.84& 0.08& 0.08& 0.19\\
    10000  & 2.21& 2.37& 3.35& 0.13& 0.11& 0.27\\
    50000  & 1.26& 1.37& 1.07& 0.62& 0.35& 1.51\\
    100000 & 0.89& 0.87& 0.67& 0.94& 1.12& 2.48
\\
    \midrule
    \multicolumn{7}{c}{\textit{Homogeneous network}} \\
    \midrule
    5000   & 3.09& 3.28& 4.45
& 0.00       & 0.00       & 0.00       \\
    10000  & 2.51& 2.69& 3.76
& 0.00       & 0.00       & 0.00       \\
    50000  & 1.40& 1.52& 0.80
& 0.00       & 0.00       & 0.00       \\
    100000 & 0.81& 0.58& 0.50
& 0.00       & 0.00       & 0.00       \\
    \bottomrule
    \end{tabular}
    \end{table}

    The numerical results demonstrate that the GP-based approach consistently outperforms the two benchmark algorithms. Compared with the coordinate descent method, GP yields systematic reductions in total system cost, with improvements ranging from approximately 1\% to over 4\%, indicating that coordinate descent often converges to suboptimal stationary points despite its stability across ten different initializations.
    
    When compared with the general-purpose nonlinear solver fmincon, GP still exhibits clear advantages in heterogeneous networks, particularly under high demand and low $\mu$. For homogeneous networks, fmincon matches the GP solutionl, primarily because the problem involves only four decision variables, making the global optimum relatively easy to identify. In addition, fmincon was initialized from ten distinct starting points. For heterogeneous networks, in the same experimental setting, the resulting solutions exhibit non-negligible variability, with differences between the minimum and maximum total system costs exceeding 10\%. In such cases, the solution with the lowest total system cost was selected for comparison with the GP benchmark. This underscores the sensitivity and limited robustness of fmincon when applied to higher-dimensional, nonconvex formulations.

    Overall, these results highlight the robustness and efficiency of the GP formulation in delivering globally optimal solutions, especially for large-scale and highly heterogeneous transit network design problems.
    
	\section{Summary and Conclusions}
	\label{Sec5}
    This study applies GP to the CA–based transit network design problem, yielding solutions with guaranteed global optimality. Numerical experiments demonstrate that conventional gradient-based algorithms and general-purpose nonlinear programming solvers cannot consistently ensure solution quality as experimental conditions vary, even when convergence is achieved from multiple initializations. In contrast, the GP-based approach provides a robust and reliable tool for large-scale transit network optimization.

    Nevertheless, the applicability of GP is constrained by its stringent requirements on posynomial objective and constraint structures, which limit its broader use in CA-based transit network design. Future research may explore extensions such as sequential geometric programming (SGP), signomial geometric programming, or mixed-integer geometric programming (MIGP), leveraging convex approximations and relaxations to expand the applicability of GP-based methods in this domain.
	
    \appendix
    \subsection{Nomenclature}\label{app1}
    In this appendix, we summarize all model variables and parameters used in the study. Table~\ref{Tab_Model_Variables_Parameters} provides the symbols of the model variables, along with their definitions, units, and corresponding values.

    \begin{table}
        \renewcommand{\arraystretch}{1.3}
        \centering
        \caption{Summary of model variables and parameters.}
        \label{Tab_Model_Variables_Parameters}
        \begin{threeparttable}
        \begin{tabular}{l p{6cm}}
            \hline
            Decision variables & Description\\
            \hline
            $\delta_i \left( x,y \right)$ &  Line density in direction $i \in I$ at location $\left( x,y \right)$ ($\frac{1}{km}$). \\
            $h_i \left( x,y \right)$  & Headway in direction $i \in I$ at location $\left( x,y \right)$ ($\frac{hr}{veh}$). \\
            \hline
            Other variables & \\
            \hline
            $q_i \left( x,y \right)$  & Vehicle flow in direction $i \in I$ at location $\left( x,y \right)$ ($\frac{veh}{hr \cdot km}$). \\
            $Z,Z_A,Z_P$  & Objective metrics: total system cost ($hr$), Agency cost ($\$$), Passenger cost ($hr$). \\
            $N_l,N_s,N_k,N_h$  & Agency metrics: line infrastructure length ($km$), number of stops,  vehicle distance traveled ($\frac{veh \cdot km}{hr_o}$),  vehicle time traveled ($\frac{veh \cdot hr}{hr_o}$). \\
            $T_a,T_w,T_r,T_t$  & Passenger metrics: access/egress cost ($hr$), wait time ($hr$), in-vehicle travel time ($hr$), transfer penalty ($hr$). \\
            \hline
            Model parameters & \\
            \hline
            $\pi_l,\pi_s,\pi_k,\pi_h$  & Unit cost of line infrastructure length ($0 \, \frac{\$}{km}$), number of stops ($0 \, \$$), vehicle-kilometers ($ 2 \, \frac{\$}{veh \cdot km}$), vehicle-hours ($40 \, \frac{\$}{veh \cdot hr}$). \\
            $ \lambda\left(x_o,y_o,x_d,y_d \right)$ & Passenger demand density with origin at $\left(x_o,y_o \right)$ and destination at $\left(x_d,y_d \right)$ ($\frac{pas}{km^4 \cdot hr_o}$). \\
            $ \lambda_{bo}^i\left(x,y \right)$ & Boarding passenger density in direction  $i \in I$ at location $\left( x,y \right)$ ($\frac{pas}{km^2 \cdot hr_o}$).\\
            $ \lambda_{al}^i\left(x,y \right)$ & Alighting passenger density in direction  $i \in I$ at location $\left( x,y \right)$  ($\frac{pas}{km^2 \cdot hr_o}$).\\
            $ \lambda_{fl}^i\left(x,y \right)$ & On-board flux density in direction  $i \in I$ at location $\left( x,y \right)$  ($\frac{pas}{km \cdot hr_o}$).\\
            $ \lambda_{tr}^i\left(x,y \right)$ & Transferring density in direction  $i \in I$ at location $\left( x,y \right)$  ($\frac{pas}{km^2 \cdot hr_o}$).\\
            $v$ & Cruising speed of transit vehicles ($25 \, \frac{km}{hr}$).\\
            $v_w$ & Walking speed of passengers ($2 \, \frac{km}{hr}$).\\
            $C$ & Capacity of a transit vehicle ($80 \, \frac{pas}{veh}$).\\
            $D$ & Total Passenger demand \\
            & ($5,000, \, 10,000, \, 50,000, \, 100,000 \, pas$). \\
            $\mu$ & Passengers' average VOT ($25, \, 20, \, 5 \,\frac{\$}{hr}$).\\
            $\beta_w$ & Perceived walking time factor of passengers (2) \\
            $\tau$ & Fixed delay per stop ($\frac{30}{3600}hr$).\\
            $\sigma$ & Transfer penalty ($\frac{60}{3600}hr$).\\
            $\Delta$ & Discretization parameter ($0.5 \,km$).\\
            \hline
        \end{tabular}
        \begin{tablenotes}
            \footnotesize
            \item $^1$ $hr$ denotes the time unit of hours, and $hr_o$ denotes the operating time of a transit system.
            \item $^2$ For each variable and parameter, the parentheses specify the unit and, where applicable, the numerical value used in the case studies.
        \end{tablenotes}
        \end{threeparttable}
    \end{table}

    \subsection{Heterogeneous demand distribution generation}\label{app2}
\begin{figure}[htbp]
    \centering
    \subfigure[Monocentric demand]{
        \includegraphics[width=0.45\linewidth]{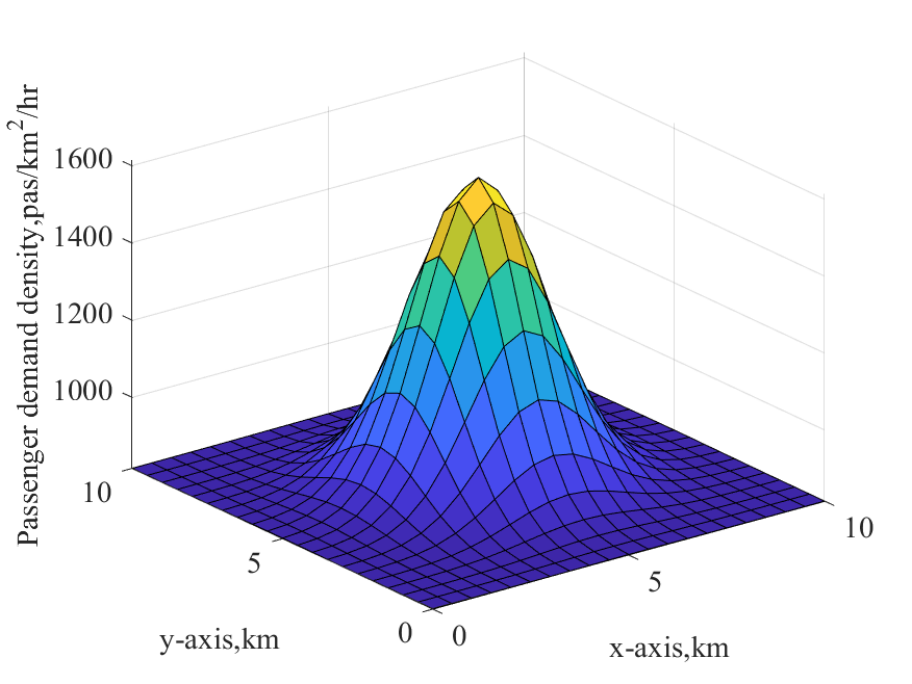}
        \label{Fig_monocentric_demand}
    }
    \subfigure[commute demand]{
        \includegraphics[width=0.45\linewidth]{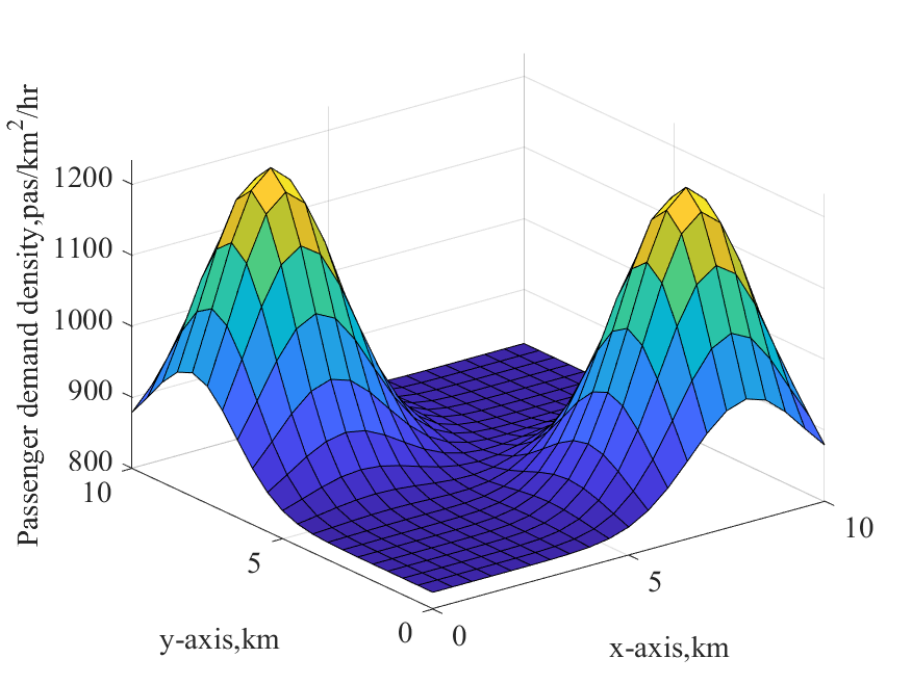}
        \label{Fig_commute_demand}
    }
    \subfigure[Chessboard1 demand]{
        \includegraphics[width=0.45\linewidth]{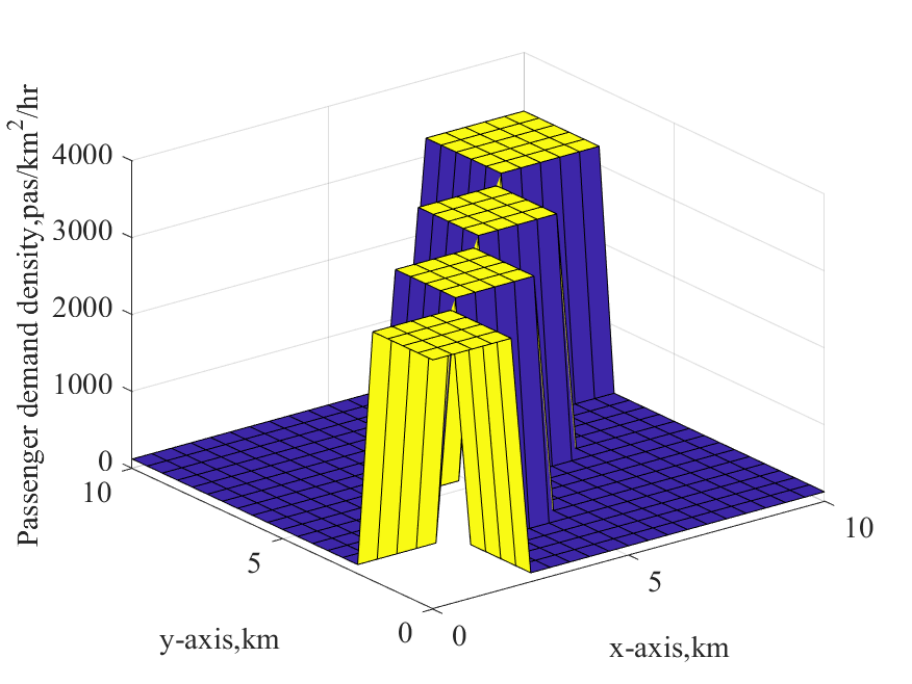}
        \label{Fig_Chessboard0_demand}
    }
    \subfigure[Chessboard2 demand]{
        \includegraphics[width=0.45\linewidth]{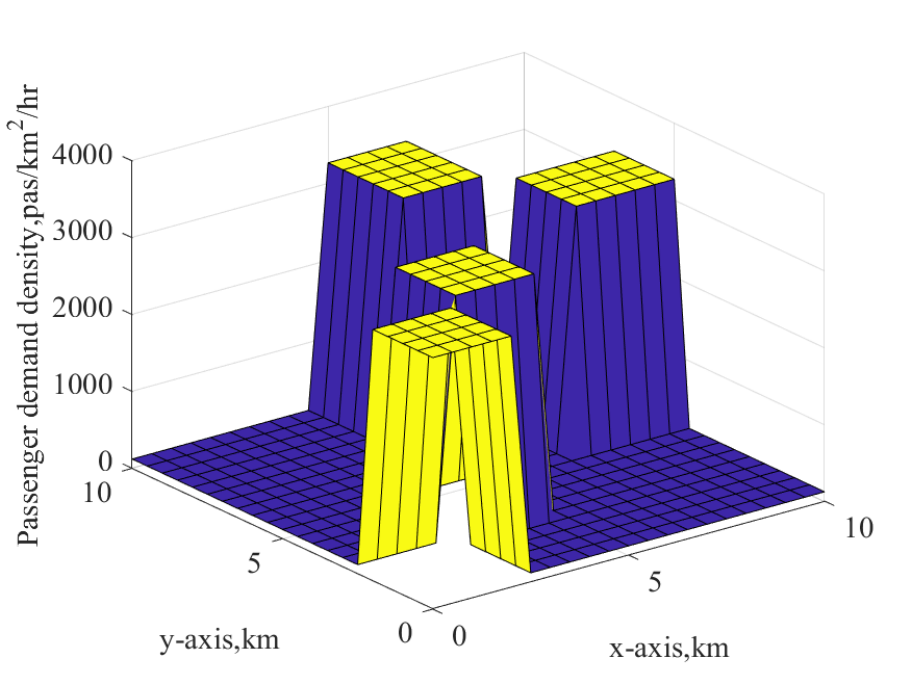}
        \label{Fig_Chessboard1_demand}
    }
    \subfigure[Chessboard3 demand]{
        \includegraphics[width=0.45\linewidth]{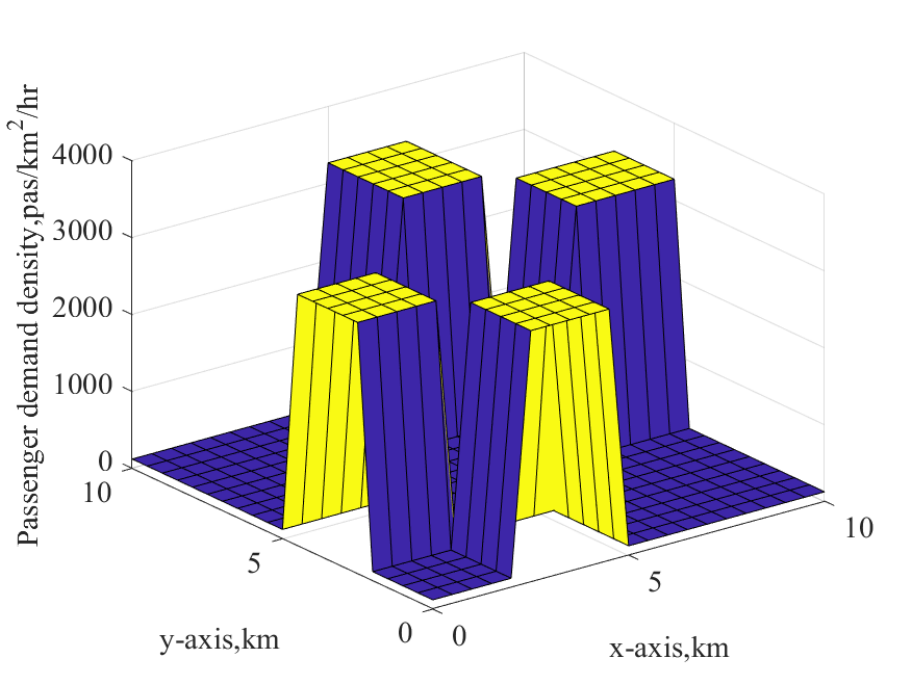}
        \label{Fig_Chessboard2_demand}
    }
    \subfigure[Chessboard4 demand]{
        \includegraphics[width=0.45\linewidth]{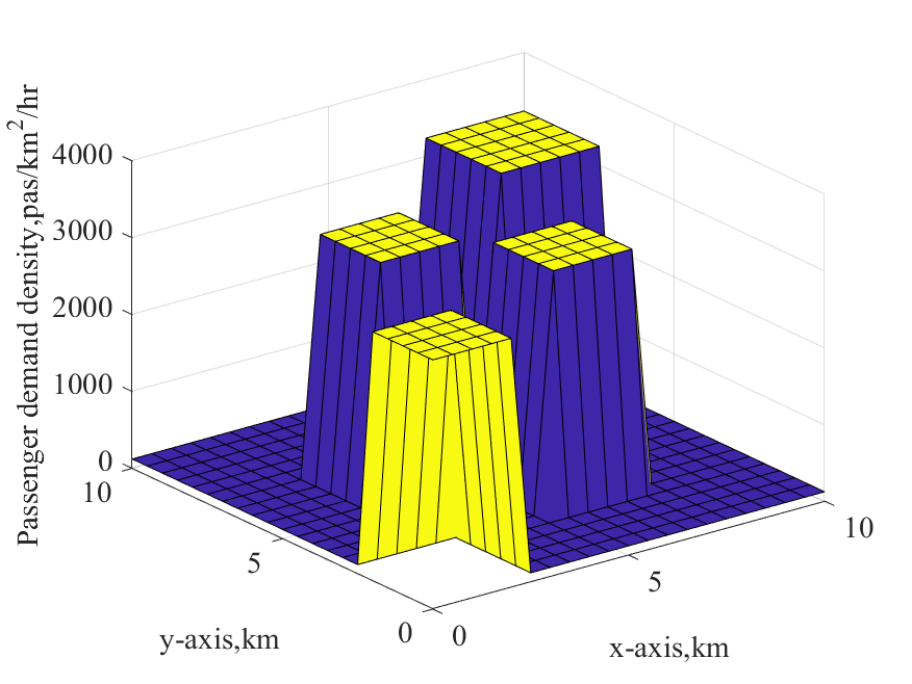}
        \label{Fig_Chessboard3_demand}
    }
    \caption{Sum of origin and destination demand densities ($D=50,000$).}
    \label{Fig_demand}
    \end{figure}

    For the monocentric and commute demand distributions, we follow the method outlined by Ouyang et al. (2014) \cite{ouyangContinuumApproximationApproach2014}, as described in Eq. (\ref{Eq_monocentric_and_commute_demand}), where $D$ represents the total passenger demand. The parameter values in Eq. (\ref{Eq_monocentric_and_commute_demand_1}) are provided in Table \ref{Tab_monocentric_and_commute_demand_density_function_parameter}. Eq. (\ref{Eq_monocentric_and_commute_demand_2}) scales the demand density function proportionally based on the given $D$. The generated results are shown in Figs. \ref{Fig_monocentric_demand} and \ref{Fig_commute_demand}.

    \begin{subequations}\label{Eq_monocentric_and_commute_demand}
    \begin{align}
    & \begin{aligned}[t]
    & \lambda'(x_o,y_o,x_d,y_d) = \\
    & \prod_{\theta=o,d}\left[ \alpha_1+\alpha_2 \sum_{\gamma=1}^{2}e^{-(\alpha_{3\gamma}x_\theta-\beta_{\theta\gamma})^2-(\alpha_{4\gamma}x_\theta-\bar{\beta}_{\theta\gamma})^2} \right] ,
    \end{aligned}
    \label{Eq_monocentric_and_commute_demand_1}\\
    & \begin{aligned}[t]
    & \lambda(x_o,y_o,x_d,y_d) = \\
    & \frac{D}{\iiiint_{0}^{|R|}\lambda'(x_o,y_o,x_d,y_d)\,dx_o\,dy_o\,dx_d\,dy_d}\lambda'(x_o,y_o,x_d,y_d) .
    \end{aligned}
    \label{Eq_monocentric_and_commute_demand_2}
    \end{align}
    \end{subequations}

    \begin{table}
    \renewcommand{\arraystretch}{1.3}
        \centering
        \caption{Parameters of the demand density function for monocentric and commute patterns.}
        \label{Tab_monocentric_and_commute_demand_density_function_parameter}
        \begin{tabular}{llllllll}
            \hline
            Demand & $\alpha_1$  & $\alpha_2$ & $\alpha_{31}$ & $\alpha_{32}$ & $\alpha_{41}$ & $\alpha_{42}$ & $\beta_{o1}$ \\
            \hline
            Monocentric & 0.0016 & 0.065 & 0.5 & 0 & 0.5 & 0 & 2.5 \\
            Commute & 0.00044 & 0.70 & 0.5 & 0 & 0.5 & 0 & 1.0 \\
            \hline
            Demand & $\beta_{o2}$  & $\beta_{d1}$ & $\beta_{d2}$ & $\bar{\beta}_{o1}$ & $\bar{\beta}_{2}$ & $\bar{\beta}_{d1}$ & $\bar{\beta}_{d2}$ \\
            \hline
            Monocentric & 0 & 2.5 & 0 & 2.5 & 0 & 2.5 & 0 \\
            Commute & 0 & 4.0 & 0 & 4.0 & 0 & 1.0 & 0 \\         
            \hline
        \end{tabular}
    \end{table}

    To further enhance the heterogeneity of passenger demand distribution, we propose a chessboard demand distribution. The core idea is to divide the network into high-density $H$ and low-density $L$ regions, which are distributed in a chessboard-like pattern across the network. The passenger demand density remains constant within each $H$ and $L$ region. The higher heterogeneity of this pattern arises from the asymmetric distribution of travel origins and destinations, which may facilitate a more realistic design. The design process for this demand distribution consists of three steps:
    
    \begin{enumerate}
    \item Identifying the locations of each $H$ and $L$ regions.
    \item Constructing passenger transfer equation set to solve for demand densities in the four transfer scenarios: high-density to high-density, $\lambda_{H \rightarrow H}$; high-density to low-density, $\lambda_{H \rightarrow L}$; low-density to high-density, $\lambda_{L \rightarrow H}$; and low-density to low-density, $\lambda_{L \rightarrow L}$. 
    
    The equation set is given in Eq. (\ref{Eq_chessboard_demand_equation_set}), where Eq. (\ref{Eq_chessboard_demand_equation_set_1}) calculates the departure and arrival demand densities, and Eq. (\ref{Eq_chessboard_demand_equation_set_2}) constrains the total demand and relevant ratios. $R_H$ and $R_L$ represent the areas of high and low demand regions, respectively. $\lambda_{bo,H}$ and $\lambda_{bo,L}$ denote the demand densities leaving the high- and low-density regions, while $\lambda_{al,H}$ and $\lambda_{al,L}$ represent the demand densities arriving in these regions. $\rho_H$ is the proportion of demand leaving and arriving in high-density areas, and $\rho_{H \rightarrow H}$ is the proportion of demand within high-density regions that corresponds to the $H \rightarrow H$ transfer scenario. In this study, both of these ratios are set to 0.9.

    \begin{subequations}\label{Eq_chessboard_demand_equation_set}
    \begin{align}
    \left\{
    \begin{aligned}
    \lambda_{bo,H} &= R_H \lambda_{H \rightarrow H} + R_L \lambda_{H \rightarrow L} \\
    \lambda_{bo,L} &= R_H \lambda_{L \rightarrow H} + R_L \lambda_{L \rightarrow L} \\
    \lambda_{al,H} &= R_H \lambda_{H \rightarrow H} + R_L \lambda_{L \rightarrow H} \\
    \lambda_{al,L} &= R_H \lambda_{H \rightarrow L} + R_L \lambda_{L \rightarrow L}
    \end{aligned}
    \right. 
    \label{Eq_chessboard_demand_equation_set_1} \\
    \left\{
    \begin{aligned}
    R_H \lambda_{bo,H} + R_L \lambda_{bo,L} &= D \\
    R_H \lambda_{al,H} + R_L \lambda_{al,L} &= D \\
    \frac{R_H \lambda_{bo,H}}{D} = \frac{R_H \lambda_{al,H}}{D} &= \rho_H \\
    \frac{R_H \lambda_{H \rightarrow H}}{\lambda_{bo,H}} &= \rho_{H \rightarrow H}
    \end{aligned}
    \right.
    \label{Eq_chessboard_demand_equation_set_2}
    \end{align}
    \end{subequations}

    The solution is given in Eq. (\ref{Eq_chessboard_demand_solution}). Notably, $\lambda_{H \rightarrow L} = \lambda_{L \rightarrow H}$, as $\frac{R_H \lambda_{bo,H}}{D} = \frac{R_H \lambda_{al,H}}{D} = \rho_H$.
    
    \begin{equation}\label{Eq_chessboard_demand_solution}
    \left\{
    \begin{aligned}
    \lambda_{H \rightarrow H} &= \frac{D \rho_H}{R_H^2(2-\rho_{H \rightarrow H})} \\
    \lambda_{H \rightarrow L} &= \lambda_{L \rightarrow H} = \frac{D \rho_H (1-\rho_{H \rightarrow H})}{R_H R_L(2-\rho_{H \rightarrow H})} \\
    \lambda_{L \rightarrow L} &= \frac{D (2-\rho_{H \rightarrow H}-3\rho_H+2\rho_H\rho_{H \rightarrow H})}{R_L^2(2-\rho_{H \rightarrow H})} \\
    \end{aligned}
    \right.
    \end{equation}
    
    \item Assigning values to $\lambda(x_o,y_o,x_d,y_d)$ based on the results from step 2) and the locations determined in step 1).
    
    \end{enumerate}
    
    Finally, we obtained four distinct chessboard demand distributions, as shown in Figs. \ref{Fig_Chessboard0_demand}, \ref{Fig_Chessboard1_demand}, \ref{Fig_Chessboard2_demand}, and \ref{Fig_Chessboard3_demand}.

    \subsection{Coordinate descent algorithm for model solution}\label{app3}

    For heterogeneous networks, the solution procedure is summarized as follow.
    
    First, initial values of decision variable $\delta_{EW}^{(0)}(y_{n'}),\,\delta_{NS}^{(0)}(x_n),\,h_{EW}^{(0)}(y_{n'}),\,h_{NS}^{(0)}(x_n)$ are randomly generated.
    
    Second, holding the line densities $\delta$ fixed,  the unconstrained optimal headways are obtained by solving the first-order optimality conditions of the total cost function. For the $E,W$ direction, the candidate headway at iteration $m$ is given by
    
    \begin{equation}\label{Eq_coordinate_descent_hetnet_h}
    \tilde{h}_{EW}^{(m)}(y_{n'})=
    \sqrt{
    \frac{
    2\delta_{EW}^{(m-1)}(y_{n'}) 
    \left( 
    \begin{aligned}
    &\textstyle \frac{\pi_k|R|}{\mu} + \frac{\pi_h|R|}{\mu v} + \\
    &\textstyle \frac{\pi_h \tau}{\mu}\Delta\sum_{n=1}^{\mathcal{N}}\delta_{NS}^{(m-1)}(x_n)
    \end{aligned}
    \right)
    }
    {\Delta\sum_{n=1}^{\mathcal{N}} \left( \lambda_{bo}^{EW}(x_n,y_{n'}) + \lambda_{al}^{EW}(x_n,y_{n'}) \right)}}
    \end{equation}
    
    Third, with the headways $h$ fixed, the optimal line densities are analogously derived from the corresponding first-order conditions. The expression for the $E,W$ direction is
    
    \begin{equation}\label{Eq_coordinate_descent_hetnet_delta}
    \begin{aligned}
    & \delta_{EW}^{(m)}(y_{n'})= \\
    &\sqrt{
    \frac{\omega_w\Delta\sum_{n=1}^{\mathcal{N}} \sum_{i'=\{EW,NS\}} \left( \lambda_{bo}^{i'}(x_n,y_{n'}) + \lambda_{al}^{i'}(x_n,y_{n'}) \right)}
    {4v_w
    \left[
    \begin{aligned}
    & \textstyle \frac{\pi_l|R|}{\mu} + \frac{\pi_s}{\mu}\Delta\sum_{n=1}^{\mathcal{N}}\delta_{NS}^{(m-1)}(x_n) + \frac{1}{\tilde{h}_{EW}^{(m)}(y_{n'})} \times \\
    & \textstyle \left( \frac{\pi_k|R|}{\mu} + \frac{\pi_h|R|}{\mu v} + \frac{\pi_h \tau}{\mu}\Delta\sum_{n=1}^{\mathcal{N}}\delta_{NS}^{(m-1)}(x_n) \right) + \\
    & \textstyle \frac{\tau\pi_h}{\mu} \Delta\sum_{n=1}^{\mathcal{N}} \left( \frac{\delta_{NS}^{(m-1)}(x_n)}{\tilde{h}_{NS}^{(m)}(x_n)} + \frac{\mu}{\pi_h} \lambda_{fl}^{NS}(x_n,y_{n'}) \right)
    \end{aligned}
    \right]
    }}
    \end{aligned}
    \end{equation}

    Fourth, to ensure that the capacity constraint is satisfied, the headways are updated accordingly. Taking the $E,W$ direction as an example, the feasible headway is obtained as

    \begin{equation}\label{Eq_coordinate_descent_hetnet_h_update}
    h_{EW}^{(m)}(y_{n'})=\min\left\{\tilde{h}_{EW}^{(m)}(y_{n'}),\frac{C \cdot \delta_{EW}^{(m)}(y_{n'})}{\lambda_{fl}^{EW}(y_{n'})}\right\}
    \end{equation}
    
    Finally, Steps 2-4 are repeated until convergence is achieved, which typically occurs within approximately ten iterations.

    For homogeneous networks, the same coordinate descent framework applies; however, due to space limitations, the detailed derivation is omitted.


    
	\bibliographystyle{IEEEtran}
	\bibliography{root} 

\end{document}